\documentclass{article}
\usepackage{amssymb,latexsym,amsmath,amsthm,mathrsfs}
\usepackage{graphicx}
\newcommand{\comment}[1]{}

\begin{document}
\title{On the values of integrals extended from the variable term $x=0$ up to
$x=\infty$\footnote{Presented to the
St. Petersburg Academy on April 30, 1781.
Originally published as
{\em De valoribus integralium a termino variabilis $x=0$ usque ad $x=\infty$ extensorum},
Institutiones calculi integralis \textbf{4} (1794), 337--345.
E675 in the Enestr{\"o}m index.
Translated from the Latin by Jordan Bell, Department of Mathematics,
University of Toronto, Toronto, Ontario, Canada. Email: jordan.bell@gmail.com}}
\author{Leonhard Euler}
\date{}
\maketitle

\S 1.\footnote{Translator: The Institutiones calculi integralis
is a collection of several of Euler's papers on the integral calculus,
and the sections of this paper are numbered from 124 to 140 in it.}
 Of those formulae which, when extended from the term $x=0$ up to
the term $x=\infty$, are assigned a finite value,
 the simplest is the circular $\int \frac{\partial x}{1+xx}$,
whose value is $\frac{\pi}{2}$,
with $\pi$ denoting the periphery for the diameter $=1$.
After this, I have also found by a singular method that
\[
\int \frac{x^{m-1}\partial x}{(1+x)^n}
\begin{bmatrix}
\textrm{from}&x=0\\
\textrm{to}&x=\infty
\end{bmatrix}
=
\frac{\pi}{n\sin\frac{m\pi}{n}}.
\]
Also by this method I have in fact obtained many other formulae of this type,
in whose differentials enter not only algebraic functions of $x$ but also
$lx$.

\S 2. Previously other formulae of this type, involving transcendental
functions, have presented themselves to me,
whose desired values seemed to refuse all methods so far known.
In particular I had searched for a curved line in which the radius of osculation
was everywhere in reciprocal proportion to the arc of the curve, so that
by putting the arc $=s$ and the radius of osculation $=r$, it would be 
$rs=aa$. 
Here it is hardly difficult to describe the figure as drawn by hand,
which gives a fair idea of what kind of figure it has.\footnote{Translator:
The paper references ``Fig. 2'', but I do not have access to it. The paper makes
sense if we assume that the figure is a clothoid (a certain spiral)
with $A$ the origin, $O$ an end point, and $C$ the abscissa of $O$, so
that
$AC$ would be the $x$ coordinate of $O$ and $CO$ would be the $y$ coordinate
of $O$.}
Doubtlessly with the initial point of the curve constituted at $A$,
from this the curve will continually bend in more, and after infinitely many spirals
will accumulate at a certain point $O$,
which may be called the {\em pole} of this curve.
My objective therefore had been to accurately investigate the location
of this pole, and for this, to seek the quantities of the coordinates
$AC$ and $CO$.

\S 3. To this end, 
having introduced the calculation of any portion $AM=s$ 
with
amplitude $=\phi$, so that it would be $r=\frac{\partial s}{\partial \phi}$
and $s\partial s=aa\partial \phi$, 
then 
\[
ss=2aa\phi, \quad \textrm{and}\quad s=a\surd 2\phi=2c\surd \phi.
\]
Hence it follows now that $\partial s=\frac{c\partial \phi}{\surd \phi}$,
from which, having put the abscissa for the arc $AP=x$ and putting
$PM=y$, it follows that
\[
x=c\int \frac{\partial \phi \cos \phi}{\surd \phi} \quad \textrm{and}\quad y=c\int \frac{\partial \phi\sin \phi}{\surd \phi}.
\]

\S 4. Then for the determination of the pole $O$,
the values of these two integral formulae are therefore required,
after they are extended from the term $\phi=0$ up to $\phi=\infty$.
Indeed I initially thought that these values could not be otherwise
obtained except by approximation, where each formula would be successively
expanded by parts; first namely from $\phi=0$ up to $\phi=\pi$;
then from
$\phi=\pi$ up to $\phi=2\pi$; next from $\phi=2\pi$ up to $\phi=3\pi$; etc.,
from which follow rapidly convergent series. However it is clear
that this operation requires long and rather tedious calculations,
which indeed I was not eager to expand. But recently I was fortunate enough
to
observe, by a singular method, that
\[
\begin{split}
&\int \frac{\partial \phi\cos\phi}{\surd \phi}
\begin{bmatrix}\textrm{from}&\phi=0\\ \textrm{to}&\phi=\infty\end{bmatrix}
=
\surd \frac{\pi}{2}, \quad \textrm{and}\\
&\int \frac{\partial \phi\sin\phi}{\surd \phi}
\begin{bmatrix}\textrm{from}&\phi=0\\ \textrm{to}&\phi=\infty\end{bmatrix}
=\surd \frac{\pi}{2};
\end{split}
\]
thus for the location of the pole $O$ that is being searched for, it would
be
\[
AC=c\surd \frac{\pi}{2} \quad \textrm{and} \quad CO=c\surd\frac{\pi}{2}.
\]

\S 5. Therefore because this method, which I have worked out here,
seems to offer no small promise, it would seem hardly disagreeable to
Geometers if I explained it here with all care. And because it can be applied
more widely than to just these formulae, I will even propose it 
in full generality. I have deduced everything just from the consideration of
this integral $\int x^{n-1}\partial x e^{-x}$, hence it is appropriate to investigate
this integral for various values of the exponent $n$.

\S 6. First of course, for the case $n=1$, for the formula
$\int \partial x e^{-x}$ the integral is clearly $1-e^{-x}$, which
vanishes in the case $x=0$
and on the other hand goes to unity by making $x=\infty$. 
Thereafter, since the differential of the formula $x^\lambda \cdot e^{-x}$ is
\[
\lambda x^{\lambda-1}\partial x \cdot e^{-x} - x^\lambda \partial x\cdot e^{-x},
\]
it will in turn be
\[
\int x^\lambda \partial x \cdot e^{-x}=\lambda \int x^{\lambda-1}\partial x\cdot e^{-x}-x^\lambda \cdot e^{-x},
\]
where the latter part vanishes both for the case $x=0$ and $x=\infty$,
but only when $\lambda>0$. Then for our terms of the integral, it will therefore
be
\[
\int x^\lambda \partial x\cdot e^{-x}=\lambda \int x^{\lambda-1}\partial x\cdot e^{-x},
\]
by means of which formula, because $\int \partial xe^{-x}=1$,
the following values of integrals may be deduced
\[
\begin{split}
&\int x\partial xe^{-x}=1\\
&\int x^2\partial x\cdot e^{-x}=1\cdot 2\\
&\int x^3\partial x\cdot e^{-x}=1\cdot 2\cdot 3\\
&\int x^4\partial x\cdot e^{-x}=1\cdot 2\cdot 3\cdot 4
\end{split}
\]
and so in general
\[
\int x^{n-1}\partial xe^{-x}=1\cdot 2\cdot 3\cdot 4\cdots (n-1),
\]
the values of which product are naturally produced when $n$ is a positive integral number;
while when $n$ is a fractional number I have previously shown how
the values can be exhibited by the quadrature of algebraic curves. It is thus
established for the case $n=\frac{1}{2}$ that its value is
$=\surd{\pi}$.

\S 7. Therefore since all the values of this infinite product
$1\cdot 2\cdot 3\cdot 4\cdots (n-1)$ can be thought of as known,
I will designate this with the letter $\Delta$, so that it would thus be
$\Delta=1\cdot 2\cdot 3\cdot 4\cdots (n-1)$,
and thus we now arrive at this remarkable integral formula
\[
\int x^{n-1}\partial x\cdot e^{-x}=\Delta,
\]
of course with the integral extended from $x=0$ to $x=\infty$;
and from this formula I have deduced everything about the case mentioned above,
where indeed some singular calculations need to be applied, 
which I will therefore carefully explain now.

\S 8. I first put $x=ky$, and because both terms of the integral
stay the same, it will be
\[
k^n \int y^{n-1}\partial y\cdot e^{-ky}=\Delta,
\]
where this formula is also extended from $y=0$ up to $y=\infty$;
then dividing this by $k^n$ we will have
\[
\int y^{n-1}\partial y\cdot e^{-ky}=\frac{\Delta}{k^n},
\]
where it should however be noted that no negative numbers can be taken for $k$,
as otherwise the formula $e^{-ky}$ would no longer vanish in the case
$y=\infty$;
and these are the only values which ought to be excluded here, so that
even imaginary values could be used in place of $k$, and then
I pursued these laborious integrations.

\S 9. Let us therefore put $k=p+q\surd{-1}$, and since it is
\begin{eqnarray*}
e^{-qy\surd -1}&=&\cos qy-\surd -1\sin qy, \quad \textrm{and}\\
e^{+qy\surd -1}&=&\cos qy+\surd -1\sin qy,
\end{eqnarray*}
our formula will now assume this form
\[
\int y^{n-1}\partial y\cdot e^{-py}(\cos qy-\surd -1 \sin qy)=
\frac{\Delta}{(p+q\surd -1)^n}.
\]
Then if we change the sign of the imaginary formulae, in a similar
way it will be
\[
\int y^{n-1}\partial y\cdot e^{-py}(\cos qy+\surd -1\sin qy)=
\frac{\Delta}{p-q\surd -1)^n}.
\]

\S 10. To help us express the values that are found in a convenient way,
let us put $p=f\cos\theta$ and $q=f\sin\theta$, and it will be
\begin{eqnarray*}
(p+q\surd -1)^n&=&f^n(\cos n\theta+\surd -1\sin n\theta)\quad \textrm{and}\\
(p-q\surd -1)^n&=&f^n(\cos n\theta-\surd -1\sin n\theta);
\end{eqnarray*}
it will be helpful to note here that $\tan\theta=\frac{q}{p}$,
whence from the assumed values $p$ and $q$ it will be $f=\surd(pp+qq)$.
In this way the first case, in the first case it would be
\[
\frac{\Delta}{(p+q\surd -1)^n}=\frac{\Delta}{f^n(\cos n\theta+\surd -1\sin n\theta)},
\]
for the second
\[
\frac{\Delta}{(p-q\surd -1)^n}=\frac{\Delta}{f^n(\cos n\theta-\surd -1\sin n\theta)}.
\]
Now if these two formulae are added, it becomes
\[
\frac{2\Delta\cos n\theta}{f^n}.
\]
On the other hand, the difference of these formulae gives
\[
\frac{2\Delta\surd -1\sin n\theta}{f^n}.
\]

\S 11. We may also add these integral formulae, and we will have
\[
\int y^{n-1}\partial y\cdot e^{-py}\cos qy=\frac{\Delta \cos n\theta}{f^n}.
\]
On the other hand, we may also subtract and divide by $2\surd -1$,
and it follows that
\[
\int y^{n-1}\partial y\cdot e^{-py}\sin qy=\frac{\Delta \sin n\theta}{f^n}.
\]
These two integral formulae are now wide open, since the numbers
$p$ and $q$ remain fully our choice, except as far as has already been observed,
that no negative number may be taken for $p$. It will thus be worthwhile
to summarize these two integral formulae in the following pair of Theorems.

\begin{center}
{\large Theorem I.}
\end{center} 

Having put $\Delta=1\cdot 2\cdot 3\cdots (n-1)$, and taking for the letters
$p$ and $q$ any positive numbers, for which one lets $\surd(pp+qq)=f$
and asks for an angle $\theta$ such that it is $\tan \theta=\frac{q}{p}$,
this remarkable integration will be obtained
\[
\int x^{n-1}\partial x\cdot e^{-px}\cos qx
\begin{bmatrix}\textrm{from}&x=0\\\textrm{to}&x=\infty\end{bmatrix}=
\frac{\Delta \cos n\theta}{f^n}.
\]

\begin{center}
{\Large Theorem II.}
\end{center}

Having put $\Delta=1\cdot 2\cdot 3\cdots (n-1)$, and taking for the letters
$p$ and $q$ any positive numbers, for which one lets $\surd(pp+qq)=f$ and asks
for an angle $\theta$ so that it is $\tan \theta=\frac{q}{p}$,
then too this remarkable integration will be obtained
\[
\int x^{n-1}\partial x\cdot e^{-px}\sin qx
\begin{bmatrix}\textrm{from}&x=0\\\textrm{to}&x=\infty\end{bmatrix}=
\frac{\Delta \sin n\theta}{f^n}.
\]

\S 12. Since for the case of the curve considered above we were
led to these integral formulae
\[
\int \frac{\partial \phi \cos \phi}{\surd \phi} \quad \textrm{and} \quad
\int \frac{\partial \phi \sin \phi}{\surd \phi},
\]
according to which it will be $n=\frac{1}{2}$, and so $\Delta=\surd \pi$,
then indeed it will be $p=0$ and $q=1$,
whence it will be $f=1$ and $\tan \theta=\frac{q}{p}=\infty$; and so $\theta=\frac{\pi}{2}$, hence
$\cos n\theta=\frac{1}{\surd 2}=\sin n\theta$. Therefore it will be
\[
\begin{split}
&\int \frac{\partial \phi \cos \phi}{\surd\phi}
\begin{bmatrix}\textrm{from}&\phi=0\\ \textrm{to}&\phi=\infty\end{bmatrix}
=\surd \frac{\pi}{2},\quad \textrm{and similarly}\\
&\int \frac{\partial \phi \sin \phi}{\surd \phi}
\begin{bmatrix}\textrm{from}&\phi=0\\ \textrm{to}&\phi=\infty\end{bmatrix}
=\surd \frac{\pi}{2}
\end{split}
\]

\S 13. It will also be worthwhile to expand the case where
$n=\frac{1}{2}$ and $\Delta=\surd \pi$ in general, and since we have
put
\[
\begin{split}
&\surd(pp+qq)=f\quad \textrm{and} \quad \frac{q}{p}=\tan \theta,\quad \textrm{it will be}\\
&\sin\theta=\frac{q}{f}\quad \textrm{and} \quad \cos\theta=\frac{p}{f}.
\end{split}
\] 
First therefore
\[
\begin{split}
&\sin\frac{1}{2}\theta=\surd \frac{1-\cos\theta}{2}=\surd \frac{f-p}{2f} \quad \textrm{and}\\
&\cos\frac{1}{2}\theta=\surd \frac{1+\cos\theta}{2}=\surd \frac{f+p}{2f};
\end{split}
\]
whence for the integral values it will be
\[
\begin{split}
&\frac{\Delta \sin\frac{1}{2}\theta}{\surd f}=\frac{\surd \pi}{f} \surd \frac{f-p}{2} \quad \textrm{and}\\
&\frac{\Delta \cos \frac{1}{2}\theta}{\surd f}=\frac{\surd \pi}{f}\cdot \surd \frac{f+p}{2}.
\end{split}
\]
Therefore we will have the following two integral formulae
\[
\begin{split}
&\int \frac{\partial x}{\surd x}e^{-px}\sin qx=\frac{\surd \pi}{f}\cdot \surd \frac{f-p}{2}\\
&\int \frac{\partial x}{\surd x}e^{-px}\cos qx=\frac{\surd \pi}{f}\cdot \surd \frac{f+p}{2}.
\end{split}
\]

\S 14. The cases in which a positive integral
number is taken for $n$, and so $\Delta$ can be completely exhibited
by integral numbers, are such that the work of reducing the integral
formulae to known quantities is straightforward, and thus such
 integrals can be exhibited in general.
This matter requires no long calculations, because our formulae
are for the simple case $x=\infty$,
yet they still are worthy of all attention. 
If however we wanted the exponent $n$ to take negative values,
these cases would demand at the start of the integration the addition
of an infinite constant, so that the integral would in fact vanish in the case 
$x=0$, and
thus the values of the integrals which we are searching for would
remain infinite, 
and thus cannot be referred to our earlier work.

\S 15. However the most remarkable case occurs here,
when $n=0$, and which demands singular skill, which we shall
therefore expand accurately.
Since we have put
\[
\Delta=1\cdot 2\cdot 3\cdot 4\cdots (n-1),
\]
let us set in a similar way
\[
\Delta'=1\cdot 2\cdot 3\cdots n, \quad \textrm{and} \quad \Delta''=1\cdot 2\cdot 3\cdots (n+1),
\]
and it will clearly be
\[
\Delta=\frac{\Delta'}{n},\quad \textrm{and}\quad \Delta'=\frac{\Delta''}{n+1},\quad \textrm{and thus} \quad \Delta=\frac{\Delta''}{n(n+1)}.
\]
Now let us take $n=\omega$,
with $\omega$ being infinitely small, and since $\Delta''=1$,
then it would be $\Delta=\frac{1}{\omega}$, whose value will thus be infinite.
Moreover, since for the first integral formula it is $\sin n\theta=\omega \theta$,
it is evident to be $\Delta \sin n\theta=\theta$; 
whence the first integral formula will be $\int \frac{\partial x}{x}e^{-px} \sin qx=\theta$, providing of course that the integral is extended from the term $x=0$ up to the term $x=\infty$.
However the value of our other integral formula $\int \frac{\partial x}{x}e^{-px}\cos qx$ will 
be infinitely large.
This case is also of great merit, so that it shall be summarized in this
singular theorem.

\begin{center}
{\Large Theorem III.}
\end{center}

\S 16. If the letters $p$ and $q$ denote any positive numbers, and then an angle
$\theta$ is asked for such that $\tan \theta=\frac{q}{p}$,
one will have the following most remarkable integration
\[
\int \frac{\partial x}{x}e^{-px} \sin qx
\begin{bmatrix}\textrm{from}&x=0\\ \textrm{to}&x=\infty\end{bmatrix}=\theta
\]
the demonstration of which theorem had without doubt only been investigated
by approximations until now.

\S 17. Moreover for the simplest case when $p=0$ and $q=1$,
which so far had seemed to surpass all techniques of calculation so far known,
because indeed in this case it would be $\tan\theta=\frac{1}{0}=\infty$, it will
be $\theta=\frac{\pi}{2}$,
from which follows the integration $\int \frac{\partial x}{x}\sin x=\frac{\pi}{2}$.
At the same time however the truth of this can not at all be doubted,
because the approximations which were done led approximately to the same
value. 
But if we compare this case to the initially considered
$\int \frac{\partial x}{\surd x}\sin x=\surd \frac{\pi}{2}$
the great similarity merits the highest attention,
since the integral of the latter is precisely the square root of the former.

\end{document}